\documentclass[12pt]{article}  %,a4paper
\usepackage{amsmath, amssymb,amscd, amsthm} 
\title{On the general solution of the Heideman-Hogan family of recurrences} 
\author{Andrew N.W.  Hone\thanks{School of Mathematics, 
Statistics and Actuarial Science, University of
Kent, Canterbury CT2 7NF, U.K. ~~E-mail: A.N.W.Hone@kent.ac.uk}$\,$ 
and Chloe Ward} 
 
\begin{document} 
 
\renewcommand{\theequation}{\arabic{section}.\arabic{equation}}

\newcommand{\beq}{\begin{equation}}  
\newcommand{\eeq}{\end{equation}}  
\newcommand{\bea}{\begin{eqnarray}}  
\newcommand{\eea}{\end{eqnarray}}  
\newcommand\la{{\lambda}}   
\newcommand\La{{\Lambda}}   
\newcommand\ka{{\kappa}}   
\newcommand\al{{\alpha}}   
\newcommand\be{{\beta}} 
\newcommand\gam{{\gamma}}     
\newcommand\om{{\omega}}  
\newcommand\tal{{\tilde{\alpha}}}  
\newcommand\tbe{{\tilde{\beta}}}   
\newcommand\tla{{\tilde{\lambda}}}  
\newcommand\tmu{{\tilde{\mu}}}  
\newcommand\si{{\sigma}}  
\newcommand\lax{{\bf L}}    
\newcommand\mma{{\bf M}}    
\newcommand\rd{{\mathrm{d}}}  
\newcommand\tJ{{\tilde{J}}}  

\newcommand\tI{{\tilde{\mathcal{I}}}}

\newcommand\SH{{\mathcal{S}}}

\newcommand\tr{{{\mathrm{tr}}\,}}

\newtheorem{thm}{Theorem}[section]

\newtheorem{propn}[thm]{Proposition}

\newtheorem{rem}[thm]{Remark}

\newtheorem{exa}[thm]{Example}

\newtheorem{cor}[thm]{Corollary}
\newtheorem{conje}[thm]{Conjecture}

\newtheorem{defn}[thm]{Definition}
\newtheorem{lem}[thm]{Lemma}

\newenvironment{prf}{\trivlist \item [\hskip 
\labelsep {\bf Proof:}]\ignorespaces}{\qed \endtrivlist} 

\newenvironment{prfco}{\trivlist \item [\hskip 
\labelsep {\bf Proof of Corollary:}]\ignorespaces}{\qed \endtrivlist}

\newcommand{\N}{{\mathbb N}}
\newcommand{\Q}{{\mathbb Q}}
\newcommand{\Z}{{\mathbb Z}}
\newcommand{\C}{{\mathbb C}}

\maketitle
 
\begin{abstract} 
We consider a family of nonlinear rational recurrences of odd order which was introduced by Heideman and Hogan. All of 
these recurrences have the Laurent property, implying that for a particular choice of initial data (all initial values set to 1) 
they generate an integer sequence. For these particular sequences, Heideman and Hogan gave a direct proof of integrality by showing that the terms of the sequence also satisfy a linear recurrence relation with constant coefficients. Here we present an analogous result for the general solution of each of these recurrences.   
\end{abstract} 

\section{Introduction}

The theory of integer sequences generated by linear recurrences has a long history in number theory, and finds many 
applications in areas such as coding and cryptography \cite{recs}, but the case of nonlinear recurrences is much less well studied. 
For some time there has been considerable interest in rational recurrence relations of the form 
\beq\label{ratrec} 
x_{n+N}\,x_n =F(x_{n+1},\ldots, x_{n+N-1}), 
\eeq 
where $F$ is a polynomial in $N-1$ variables, 
which surprisingly generate integer sequences. Several   quadratic recurrences of this kind were found by Somos, and 
this inspired others to find new examples, as described in the articles by Gale \cite{gale}. An important early observation was that 
if (\ref{ratrec}) has the Laurent property, meaning that it generates Laurent polynomials in the initial values with 
integer coefficients, i.e. 
$$
x_n\in\Z [x_0^{\pm 1},x_1^{\pm 1},\ldots,x_{N-1}^{\pm 1}]\qquad \forall n \in \Z ,
$$ 
then an integer sequence is generated automatically by choosing the initial values to be $x_0=x_1=\ldots=x_{N-1}=1$. 
Subsequently, as an offshoot of their development of cluster algebras, Fomin and Zelevinsky introduced the Caterpillar Lemma \cite{fz02}, which is a useful tool  for proving the Laurent property for many recurrences of the form (\ref{ratrec}). 
In the special case where $F$ is a sum of two monomials, Fordy and Marsh explained how such recurrences arise from cluster 
algebras associated with quivers that are periodic under cluster mutations \cite{fordymarsh09}, while for more general $F$ 
a range of examples were found recently by Alman et al. \cite{alman}, who considered mutation periodicity  in the broader context of Laurent phenomenon (LP) algebras \cite{lp}.

In this paper we are concerned with a particular family of nonlinear recurrences of odd order 
$N=2k+1$, given by 
\beq\label{heiho} 
x_{n+2k+1}\, x_n = x_{n+2k}\, x_{n+1}+ a(x_{n+k}+x_{n+k+1}), 
\eeq 
where $a$ is a non-zero parameter. This family was introduced in the case $a=1$ 
 by Heideman and Hogan \cite{heiho}, who proved that the sequence generated by (\ref{heiho}) with 
the initial values $x_0=x_1=\ldots = x_{2k}=1$ consists entirely of integers.  
%sequence results. 
(By rescaling $x_n \to a x_n$, the  parameter $a$ can always
be removed, but it will be useful to retain it here for bookkeeping purposes.)
One way to see the integrality of this particular sequence  is to show that (\ref{heiho}) has the 
Laurent property, which was noted in \cite{heiho} and proved in \cite{hogan}; more recently, the family (\ref{heiho}) 
was rediscovered in a search for period 1 seeds in LP algebras - see Theorem 3.10 in \cite{alman}.
However, Heideman and Hogan's original proof of integrality was based on the following result.
\begin{thm}\label{mainhh}
The terms of the sequence generated by the recurrence (\ref{heiho})  
with initial values $x_j=1$ for $j=0,1,\ldots, 2k$ and $a=1$ 
satisfy the linear relation 
\beq\label{klin}
x_{n+6k}-(2k^2+8k+4)\, (x_{n+4k}-x_{n+2k}) - x_n=0
\eeq 
for all $n\in\Z$. 
\end{thm} 
The integrality result in \cite{heiho} is proved by starting from $x_j=1$, $0\leq j \leq 2k$, then 
 determining the explicit form of  the next $4k$ values $x_j$, 
$2k+1\leq j\leq 6k$ obtained by iterating  (\ref{heiho}) for $0\leq n \leq 4k-1$, 
where the value of $x_{6k}$ is used to verify that (\ref{klin}) holds for $n=0$, 
and finally showing by induction that if (\ref{klin}) is assumed to hold for all $n\geq 0$ then (\ref{heiho}) 
also holds for all $n\geq 4k$. This particular sequence is also symmetric under reversal, in the sense 
that 
\beq\label{rsym} 
x_{-n}=x_{n+2k}  \qquad \forall n\in\Z .
\eeq 
 In that case, the efficacy of this inductive approach can be seen directly from an operator identity connecting the 
linear operator in (\ref{klin}) with the nonlinear equation (\ref{heiho}), which can be rewritten in the form 
$\xi_n=0$, where 
\beq\label{xidet} 
\xi_n:=\left|\begin{array}{cc} x_n & x_{n+2k} \\ 
x_{n+1} &  x_{n+2k+1} 
\end{array}\right| - a (x_{n+k}+x_{n+k+1}). 
\eeq 
\begin{lem}\label{mlem}
Let $\SH$ denote the shift operator, such that $\SH\, x_n=x_{n+1}$ for all $n$, and let
$$ 
{\cal L}=\SH^{6k}-K(\SH^{4k}-\SH^{2k})-1, 
$$ 
where  $K$ is some fixed constant. Then 
\beq\label{opid} 
{\cal L} \, \xi_n = {\cal M}_n \cdot {\cal L} \, x_n, 
\eeq 
where $ {\cal M}_n$ is the linear operator 
$$ 
{\cal M}_n = x_{n+6k}\, \SH^{2k+1}-x_{n+6k+1}\,\SH^{2k}-x_{n+2k}\,\SH+x_{n+2k+1}-a(\SH^{k+1}+\SH^k).
$$ 
\end{lem} 

The main result of this paper is the analogue of Theorem \ref{mainhh} for the case of arbitrary initial data. 
\begin{thm}\label{main}
The iterates of the recurrence (\ref{heiho}) satisfy the linear relation 
\beq\label{gklin}
x_{n+6k}-K\, (x_{n+4k}-x_{n+2k}) - x_n=0
\eeq 
for all  $n\in\Z$, where 
\beq\label{kinv} 
K =P^{(0)}+a\, P^{(1)}+a^2\, P^{(2)}, 
\eeq  
with 
$$ 
P^{(0)}=1+\frac{x_0}{x_{2k}}+\frac{x_{2k}}{x_0}, 
$$ 
$$ 
P^{(1)}=%\left(
\Big(1+\frac{x_{2k}}{x_0}\Big)\sum_{j=1}^k\frac{x_{j-1}+x_j}{x_{j+k-1}x_{j+k}}
+ 
\Big(1+\frac{x_0}{x_{2k}}\Big)\sum_{j=1}^k\frac{x_{j+k-1}+x_{j+k}}{x_{j-1}x_j} %\right)  
, 
$$ 
$$ 
%\begin{array}{rcl} 
P^{(2)} = \frac{1}{x_kx_{2k}}+\sum_{j=0}^{k-1}\frac{1}{x_j}\left(\frac{1}{x_{j+k}}+\frac{1}{x_{j+k+1}}\right)  
%\\ 
 + \sum_{\ell=1}^{k-1}\sum_{m=1}^\ell\frac{(x_\ell+x_{\ell+1})(x_{k+m-1}+x_{k+m})}{x_{k+\ell}x_{k+\ell+1}x_{m-1}x_m}.
$$ 
\end{thm} 

In principle, it is possible to 
prove the above result directly by using the identity (\ref{opid}) in 
Lemma \ref{mlem} and adapting the argument from \cite{heiho}. 
To do so one should take $2k+1$  initial values $x_0,\ldots,x_{2k}$ for the nonlinear recurrence (\ref{heiho}), 
require the $4k$ vanishing conditions $\xi_0=\xi_1=\ldots=\xi_{4k-1}=0$ which fix 
$6k$ initial values $x_0,x_1,\ldots,x_{6k-1}$ for the linear 
equation ${\cal L} x_n=0$ together with the value of $K$, determined 
as 
\beq\label{kval} 
K=\frac{x_{6k}-x_0}{x_{4k}-x_{2k}}, 
\eeq 
and then further verify that $\xi_j=0$ for a total of $6k$ adjacent values of $j$ (including the range $0\leq j\leq 4k-1$ 
already assumed); this implies that  the corresponding  solution of the initial value problem for ${\cal L} \xi_n=0$ is the 
zero solution $\xi_n=0$ for all $n$. Heideman and Hogan made this argument  effective with the use of computer algebra, 
which they used  (for arbitrary $k$) to calculate explicit expressions for the values of $x_{2k+1},...,x_{6k}$ corresponding 
to $x_0=x_1=\ldots=x_{2k}=1$, and hence to verify the value $K=2k^2+8k+4$ in (\ref{klin}) 
and other necessary identities; 
they also implicitly used the fact that this special sequence has the reversal symmetry (\ref{rsym}) (although this fact was not stated in \cite{heiho}), which means that once $\xi_j=0$ holds for $0\leq j\leq 4k-1$ it automatically holds for  
$-2k+1\leq j\leq -1$ as well, so it is enough to verify in addition that $\xi_{4k}=0$ in order to show that $\xi_n=0$ for all $n$ 
by induction. However, this argument is much harder to apply to the case of generic initial data, because the 
corresponding sequence need not have the symmetry (\ref{rsym}), so here we prefer to  adopt a different approach. 
Nevertheless, we are able to exploit the fact that the 
recurrence (\ref{heiho}) is  itself reversible in the sense of \cite{lr}, making the proof below much simpler than it might be  otherwise. 

The result (\ref{main}) can be restated as saying that the recurrence (\ref{heiho}) is linearizable, with the coefficient 
$K$ appearing in the linear relation (\ref{gklin}) being a conserved quantity (this terminology is explained in more detail in the next section); the general solution for the case $k=1$ was already covered in \cite{honechapt}.  There are many other examples of nonlinear recurrences that are linearizable, which arise in diverse contexts ranging from cluster algebras associated with affine A-type Dynkin diagrams \cite{fordymarsh09, fordy10, fordyhone12, keller}, to 
frieze relations \cite{assem}, Q-systems for 
characters in representation theory \cite{difk}, and period 1 seeds in LP algebras \cite{alman, honeward}. In all these examples, the 
key to obtaining the linear recurrences is provided by certain determinantal identities for discrete Wronskians. The fact 
that  (\ref{heiho}) can be rewritten as the vanishing of the expression (\ref{xidet}) involving a $2\times 2$ determinant 
permits a linear relation to be derived in a straightforward way, although it turns out that this approach is insufficient to 
obtain the precise form of (\ref{gklin}).

In the next section we provide the proof of Theorem \ref{main}, and in section 3 we present various corollaries, before making some conclusions. 

\section{Proof of the main theorem} 

\setcounter{equation}{0}

Before proceeding with the proof, we give some discussion of terminology, and describe properties of (\ref{heiho}) that will be useful later on.  First of all, note that iterating the nonlinear recurrence is equivalent to iterating a birational map in 
dimension $2k+1$, namely 
\beq\label{phi} 
\varphi: \qquad (x_0,x_1,\ldots,x_{2k}) \mapsto \left(x_1,x_2,\ldots,\frac{x_1x_{2k}+a(x_{k+1}+x_k)}{x_0}\right). 
\eeq
If we always use this map to iterate then it is useful to regard the terms in the sequence  $(x_n)_{n\in\Z}$ 
as rational functions 
(in fact, Laurent polynomials, but we will not need this) in the initial coordinates $x_0,x_1,\ldots,x_{2k}$ and $a$, 
obtained by the pullback of  $\varphi$ (or its inverse) applied to these variables, so that 
$$ 
(\varphi^*)^n x_0 = x_n \qquad \forall n\in\Z, 
$$
with $(\varphi^{-1})^*=(\varphi^*)^{-1}$. We say that a non-constant function $F(x_0,x_1,\ldots,x_{2k})$ is 
a conserved quantity, or first integral, for the map $\varphi$ if it is invariant under pullback, i.e. $\varphi^*F=F\cdot \varphi=F$, 
and we say that it is a $p$-invariant if it is periodic with period $p$, i.e. $(\varphi^*)^pF=F$.

From Theorem  3.10 in \cite{alman}, the map 
can be factored as $\varphi = \rho\cdot \mu$, where $\rho$ is a cyclic permutation of the coordinates and $\mu$ is 
a mutation in an LP algebra, but more interesting for our purposes is the fact that it is a reversible map (it has the discrete 
analogue of time-reversal symmetry \cite{lr}), meaning that it is conjugate to its own inverse. 
\begin{lem}\label{rev} 
The map $\varphi$ satisfies 
$$ \varphi= 
\si\cdot \varphi^{-1}\cdot \si, 
$$ 
where the reversing symmetry $\si$ is the involution 
$$ \si: \qquad (x_0,x_1,\ldots,x_{2k}) \mapsto (x_{2k},x_{2k-1},\ldots,x_0).$$ 
\end{lem}
\noindent 
Reversibility means that the reversing symmetry $\si$  can extended to the level of the whole sequence $(x_n)$ by pullback, 
so that it acts according to 
\beq\label{revsym} 
\si^* x_n = x_{2k-n} \qquad \forall n \in\Z. 
\eeq

In order to obtain linear relations for the terms of the sequence, it will be convenient to introduce the 
$3\times 3$  discrete Wronskian matrix 
\beq\label{psi}
\Psi_n:=\left(\begin{array}{ccc} x_n & x_{n+2k} & x_{n+4k} \\
x_{n+1} & x_{n+2k+1} & x_{n+4k+1} \\
x_{n+2} & x_{n+2k+2} & x_{n+4k+2} 
\end{array}\right), 
\eeq
which has $2\times 2$ minors of the form appearing in (\ref{xidet}).
\begin{propn}\label{det} 
The determinant 
$$ 
\delta_n:=\det \Psi_n 
$$ 
%of the $3\times 3$  matrix (\ref{psi}) 
is a $k$-invariant for the map $\varphi$. 
\end{propn} 
\begin{prf} 
Using Dodgson condensation \cite{dodgson} (also known as the Desnanot-Jacobi identity) to expand the $3\times3$ determinant in terms of its 
$2\times 2$ connected minors yields 
$$\begin{array}{rcl} 
x_{n+2k+1}\,\delta_n &  =&  \left|\begin{array}{cc} \xi_n +as_{n+k} &  \xi_{n+2k} +as_{n+3k} \\ 
 \xi_{n+1} +as_{n+k+1} &  \xi_{n+2k+1} +as_{n+3k+1}
\end{array}\right|  
\\
& = & L_n +a^2  \left|\begin{array}{cc} s_{n+k} &  s_{n+3k} \\ 
 s_{n+k+1} &  s_{n+3k+1}
\end{array}\right| , 
\end{array} 
$$
where $s_n=x_n+x_{n+1}$, and the quantity $L_n$ is a sum of homogeneous linear and quadratic terms in 
$\xi_j$ for certain $j$. A direct calculation then shows that 
$$ 
x_{n+2k+1}\, x_{n+3k+1}(\delta_{n+k}-\delta_n) =   x_{n+2k+1}L_{n+k} -x_{n+3k+1}L_n + \Delta_n,  
$$ 
where 
$$ 
 \Delta_n 
=  a^2\Big(s_{n+2k+1}\,\xi_{n+2k}+s_{n+2k}\,\xi_{n+2k+1}-s_{n+3k+1}\,\xi_{n+k}-s_{n+3k}\,\xi_{n+k+1}\Big), 
$$ 
which clearly vanishes, along with $L_n$ and $L_{n+k}$, if $\xi_j=0$ for all $j$. Therefore $\delta_{n+k}=(\varphi^*)^k\delta_n = \delta_n$ for all 
$n$, as required. 
\end{prf} 
\begin{rem} \label{van} 
Working in the ambient field of rational functions, that is 
$\C(x_0,x_1,\ldots,x_{2k},a)$, and using explicit expressions for the first few iterates (see  below) 
it can be verified directly that $\delta_0$ and all its shifts $\delta_1,\ldots,\delta_{k-1}$ are non-zero 
rational functions (actually, Laurent polynomials), e.g. 
$$ 
\delta_0=\delta_{-2k} = 
\left|\begin{array}{ccc} x_{-2k} & x_{0} & x_{2k} \\
x_{-2k+1} & x_{1} & x_{2k+1} \\
x_{-2k+2} & x_{2} & x_{2k+2} 
\end{array}\right| 
$$ 
can be calculated from the formulae in  Lemma \ref{explicit}, and by periodicity none of the shifts 
$(\varphi^*)^n\delta_0$ can be identically zero (as a rational function). %the parameter $a$
\end{rem}
\begin{cor} \label{lins}
The  determinant  of the $4\times 4$ discrete Wronskian matrix 
$$ 
\hat{\Psi}_n :=\left(\begin{array}{cccc} x_n & x_{n+2k} & x_{n+4k} & x_{n+6k} \\
x_{n+1} & x_{n+2k+1} & x_{n+4k+1}& x_{n+6k+1} \\
x_{n+2} & x_{n+2k+2} & x_{n+4k+2} & x_{n+6k+2}\\ 
x_{n+3} & x_{n+2k+3} & x_{n+4k+3} & x_{n+6k+3}\\ 
\end{array}\right) 
$$ 
is zero. 
\end{cor} 
\begin{prfco} 
Using Dodgson condensation once more to expand the $4\times4$ determinant in terms of its $3\times 3$ connected minors, 
which are shifts of the determinant of  (\ref{psi}), gives 
$$ 
\det\hat{\psi}_n \, \left| \begin{array}{cc} 
x_{n+2k+1} & x_{n+4k+1} \\
x_{n+2k+2} & x_{n+4k+2} 
\end{array}\right| 
=  \left| \begin{array}{cc} 
\delta_{n} & \delta_{n+2k} \\
\delta_{n+1} & \delta_{n+2k+1} 
\end{array}\right| 
=  \left| \begin{array}{cc} 
\delta_{n} & \delta_{n} \\
\delta_{n+1} & \delta_{n+1} 
\end{array}\right| =
0
$$ by Proposition \ref{det}. 
\end{prfco} 
The above results are almost, but not quite, sufficient to produce the linear relation in Theorem \ref{main}. 
\begin{thm}\label{lnrs} The iterates of the nonlinear recurrence (\ref{heiho}) satisfy the linear recurrence 
\beq\label{k12} x_{n+6k}-K^{(1)}\, x_{n+4k}+K^{(2)}\,x_{n+2k} - x_n=0, 
\eeq 
where $K^{(1)},K^{(2)}$ are conserved quantities with 
\beq\label{k12rev} 
K^{(2)}=\sigma^* K^{(1)}, 
\eeq
as well as the linear recurrence 
 \beq\label{abg} x_{n+3}-\gam_n\, x_{n+2}+\be_n\, x_{n+1} - \al_n \,x_n=0, 
\eeq
 where $\al_n$ is a $k$-invariant and $\be_n,\gam_n$ are $2k$-invariants.
\end{thm} 
\begin{prf} An element of the kernel of $\hat{\Psi}_n$ is given by a column vector ${\bf v}_n=(-K^{(3)},K^{(2)},-K^{(1)},1)^T$, where the  
last entry has been scaled to 1 (which is valid since $\delta_n$ is non-vanishing by Remark \ref{van}), and a priori the other entries $K^{(j)}$ depend on $n$. The first three rows of the equation $\hat\Psi_n{\bf v}_n=\mathbf {0}$ give a linear system for the $K^{(j)}$, and by Cramer's rule the solution is 
$$%\beq\label{ksols} 
K^{(1)}=\frac{1}{\delta_n} \left|\begin{array}{ccc} x_n & x_{n+2k} & x_{n+6k} \\
x_{n+1} & x_{n+2k+1} & x_{n+6k+1} \\
x_{n+2} & x_{n+2k+2} & x_{n+6k+2} 
\end{array}\right|,  \,
K^{(2)}=\frac{1}{\delta_n}  \left|\begin{array}{ccc} x_n &  x_{n+4k} & x_{n+6k} \\
x_{n+1} &  x_{n+4k+1} & x_{n+6k+1} \\
x_{n+2} &  x_{n+4k+2} & x_{n+6k+2} 
\end{array}\right| 
$$%\eeq  
and $ K^{(3)}=\delta_n^{-1}\delta_{n+2k}=1$.
The last three rows of  $\hat\Psi_n{\bf v}_n=\mathbf {0}$ give the same linear system with all indices shifted by 1, implying that $K^{(1)}$ and $K^{(2)}$ are independent of $n$. Now applying $\si^*$ to (\ref{k12}), replacing 
$n\to -n-4k$ and adding the result back to the original relation produces
$$
(\si^*K^{(2)} -K^{(1)})\, x_{n+4k}+(K^{(2)}-\si^*K^{(1)})\,x_{n+2k} =0
$$
for all $n$, hence 
(\ref{k12rev}) must hold. The same argument applied to 
the kernel of  the transpose matrix $\hat{\Psi}_n^T$ yields the relation 
(\ref{abg}) where
$$%\beq\label{ksols} 
\be_n=\frac{1}{\delta_n} \left|\begin{array}{ccc} x_n & x_{n+2} & x_{n+3} \\
x_{n+2k} & x_{n+2k+2} & x_{n+2k+3} \\
x_{n+4k} & x_{n+4k+2} & x_{n+4k+3} 
\end{array}\right|,  \,
\gam_n=\frac{1}{\delta_n}  \left|\begin{array}{ccc} x_n &  x_{n+1} & x_{n+3} \\
x_{n+2k} &  x_{n+2k+1} & x_{n+2k+3} \\
x_{n+4k} &  x_{n+4k+1} & x_{n+4k+3} 
\end{array}\right| 
$$%\eeq 
are $2k$-invariants and 
$
\al_n=\delta_n^{-1}\delta_{n+1}
$
is a $k$-invariant. \end{prf} 
By considering the monodromy of the linear equation (\ref{abg}) with periodic coefficients, the coefficients in (\ref{k12}) for 
can be written in terms of $\al_j,\be_j, \gam_j$. In terms of the matrix $\Psi_n$, the system (\ref{abg}) implies that 
$$ 
\Psi_{n+1}={\bf L}_n\, \Psi_n, \quad {\bf L}_n = \left(\begin{array}{ccc}0 & 1 & 0 \\ 0 & 0 & 1 \\ 
\al_n & -\be_n & \gam_n 
\end{array}\right), 
 \quad {\bf L}_n^{-1} = \frac{1}{\al_n}\left(\begin{array}{ccc}\be_n & -\gam_n & 1 \\ \al_n & 0 & 0 \\ 
0 & \al_n & 0
\end{array}\right), 
$$ 
so that 
$$ 
\Psi_{n+2k}={\bf M}_n\, \Psi_n, \quad {\bf M}_n={\bf L}_{n+2k-1} \, {\bf L}_{n+2k-2}\cdots {\bf L}_{n},
$$ 
while on the other hand 
$$
\Psi_{n+2k}=\Psi_n\, \tilde{{\bf L}}, \qquad \tilde{{\bf L}} = \left(\begin{array}{ccc}0 & 0 & 1 \\ 1 & 0 & -K^{(2)} \\ 
0 & 1 & K^{(1)} 
\end{array}\right), 
 \quad \tilde{{\bf L}}^{-1} = \left(\begin{array}{ccc}K^{(2)} & 1 & 0 \\ - K^{(1)}  & 0 & 1 \\ 
1 & 0 & 0
\end{array}\right). 
$$
Then noting that $\tr \tilde{{\bf L}}=\tr {\bf M}_n$ and $\tr \tilde{{\bf L}}^{-1}=\tr {\bf M}_n^{-1}$, together with 
the observation that $\prod_{j=1}^k \al_j=1$, we have the following.
\begin{propn}\label{monod}
The conserved quantities $ K^{(1)}$ and $ K^{(2)}$ can be written as polynomials in  $\al_j,\be_j, \gam_j$, given by 
$$ 
 K^{(1)}=\tr {\bf L}_{n+2k-1} \, {\bf L}_{n+2k-2}\cdots {\bf L}_{n}, \qquad 
K^{(2)}=\tr {\bf L}_{n}^{-1} \, {\bf L}_{n+1}^{-1}\cdots {\bf L}_{n+2k-1}^{-1}. 
$$
\end{propn}
According to the result we are aiming for, Theorem \ref{main}, we expect to find $K^{(1)}=K^{(2)}=K$, a Laurent 
polynomial in $x_0,x_1,\ldots,x_{2k}$, which by (\ref{k12rev}) must be invariant under the reversal symmetry. However, 
none of the formulae for  $K^{(1)}, K^{(2)}$ obtained so far make this coincidence manifest, and none of them immediately yield a Laurent polynomial. Indeed, the preceding result is somewhat mysterious, since direct calculations for small values of $k$ reveal that   $\al_j,\be_j, \gam_j$ are not  Laurent polynomials themselves. In order to prove the main result, we calculate explicit formulae for $2k$ iterates on either side of the initial data, which allows  us to obtain $K$ as a quadratic 
polynomial in $a$,  by using (\ref{kval}) with suitably shifted indices. 
\begin{lem}\label{explicit} 
The first $2k$ terms obtained by iterating (\ref{heiho}) forwards from the initial values 
$x_0,x_1,\ldots,x_{2k}$ are given by 
$$ \begin{array}{rclcl} 
x_{2k+j}& =& x_0^{-1}x_j \, x_{2k}&+ & a \, F_{2k+j}^{(1)}, \\ 
x_{3k+j}& =& x_0^{-1}x_{k+j}\, x_{2k}&+ &a \, F_{3k+j}^{(1)}+ a^2\,F_{3k+j}^{(2)} , \quad 1\leq j\leq k, 
\end{array} 
$$  
where the coefficients of the linear and quadratic terms in $a$ are specified by 
$$  \begin{array}{rcl}
F_{2k+j}^{(1)}& = & x_j\sum_{\ell=1}^j (x_{\ell-1}x_\ell)^{-1} (x_{k+\ell-1}+x_{k+\ell}), \\ 
F_{3k+j}^{(1)}& = & x_0^{-1}x_{k+j}\, x_{2k} \sum_{\ell=1}^j (x_{k+\ell-1}x_{k+\ell})^{-1} (x_{\ell-1}+x_{\ell}) 
+x_k^{-1}x_{k+j}F_{3k}^{(1)}, \\ 
F_{3k+j}^{(2)}& = & x_{k+j} \sum_{\ell=1}^j  (x_{k+\ell-1}x_{k+\ell})^{-1}
\Big(F_{2k+\ell-1}^{(1)}+F_{2k+\ell}^{(1)}\Big),
\end{array} 
$$ 
for the same range of the index $j$, with $F_{2k}^{(1)}=0$. 
The  first $2k$ terms obtained by iterating (\ref{heiho}) backwards   from the same  initial values are 
$$ \begin{array}{rclcl} 
x_{-j}& =& x_{2k}^{-1}x_{2k-j} \, x_{0}&+ & a \, F_{-j}^{(1)}, \\ 
x_{-k-j}& =& x_{2k}^{-1}x_{k-j}\, x_{0}&+ &a \, F_{-k-j}^{(1)}+ a^2\,F_{-k-j}^{(2)} , \quad 1\leq j\leq k, 
\end{array} 
$$  
where, for the same range of $j$ values,  
$$ 
 F_{-j}^{(1)}= \si^* F_{2k+j}^{(1)}, \qquad F_{-k-j}^{(1)}=\si^* F_{3k+j}^{(1)}, \qquad F_{-k-j}^{(2)}= \si^* F_{3k+j}^{(2)}.
$$
\end{lem} 
\begin{prf} The first $2k+1$ iterations of (\ref{heiho}), either forwards or backwards, only require multiplication and addition of previous terms, as well as division by one of $x_0,x_1,\ldots,x_{2k}$, so for the division there is no need to consider any cancellations between numerator and denominator (which are required for the Laurent property to hold at subsequent steps). It is plain to see that the first $k$ terms produced by iterating are linear in $a$, while the next $k$ terms are quadratic. Expanding $x_{2k+j}$ in $a$ and 
substituting into (\ref{heiho}) with $n=j-1$ it is 
straightforward to obtain the leading order term recursively, while the coefficient of the linear term satisfies the recursion 
$$ 
x_j^{-1}F_{2k+j}^{(1)}-x_{j-1}^{-1}F_{2k+j-1}^{(1)} = (x_{j-1}x_j)^{-1}(x_{k+j-1}+x_{k+j}), \quad 1\leq j \leq k, 
$$
which can be summed telescopically (with $F_{2k}^{(1)}=0$) to obtain the above formula for $F_{2k+j}^{(1)}$.
Similarly, expanding in $a$ for the next $k$ iterations, the leading order term is found immediately, while for the term 
linear in $a$ the recursion is 
$$ 
x_{k+j}^{-1}F_{3k+j}^{(1)}-x_{k+j-1}^{-1}F_{3k+j-1}^{(1)} = (x_{k+j-1}x_{k+j}x_0)^{-1} x_{2k}
(x_{j-1}+x_{j}),  \quad 1\leq j \leq k. 
$$
which immediately yields the above expression for $F_{3k+j}^{(1)}$; and for the quadratic term the formula for the 
coefficient is found by solving the recursion 
$$ 
x_{k+j}^{-1}F_{3k+j}^{(2)}-x_{k+j-1}^{-1}F_{3k+j-1}^{(2)} = (x_{k+j-1}x_{k+j})^{-1}
\Big(F_{2k+j-1}^{(1)}+F_{2k+j}^{(1)}\Big), \quad 1\leq j \leq k. 
$$
Similarly,  iterating $2k$ steps backwards  produces the images of the forward iterates under the action of the reversing map, so the 
formulae for $x_{-j}$ and $x_{-k-j}$ follow by direct application of $\si^*$.
\end{prf} 
\begin{propn}\label{kexpr} 
If $K$ is defined by
\beq\label{kform}
K = \frac{x_{4k}-x_{-2k}}{x_{2k}-x_0} ,
\eeq 
where $x_{4k}$ and $x_{-2k}$ are given as Laurent polynomials in $x_0,x_1,\ldots,x_{2k}$ and $a$ according to 
Lemma \ref{explicit}, then it is given explicitly by  (\ref{kinv}).
\end{propn}
\begin{prf} 
The formula (\ref{kinv}) is readily checked at each order in $a$. At leading order this is trivial, while at order $a$ and $a^2$ 
the identities 
$$ 
(x_{2k}-x_0) \, P^{(j)} = F_{4k}^{(j)} -  F_{-2k}^{(j)} 
$$ 
are seen to hold for $j=1,2$.
\end{prf}
\begin{thm}\label{fin} The quantity $K$ in (\ref{kinv}) is a first integral for the map $\varphi$. 
\end{thm}
\begin{prf}
To verify that $K$ is a conserved quantity, observe that 
$$ 
\varphi^*K = \frac{x_{4k+1}-x_{-2k+1}}{x_{2k+1}-x_1}
$$ 
from (\ref{kform}), and this is equal to $K$ if and only if ${\cal L}\, x_{-2k+1}=0$,
where ${\cal L}$ is the operator in Lemma \ref{mlem}. Now ${\cal L}\, x_{-2k}$ vanishes by (\ref{kform}),
so  from the identity 
$$ \begin{array}{rcl}
x_{2k}{\cal L}x_{-2k+1}-x_{2k+1}  {\cal L} x_{-2k} -\xi_{2k}+K\xi_0 & = &
a(x_{3k}+x_{3k+1})-Ka(x_{k}+x_{k+1}) 
\\ && +\,x_{2k+1}x_{-2k}-x_{2k}x_{-2k+1}
\end{array}
$$ 
we see that it is sufficient to check that the right-hand above is zero. Substitution of the explicit expressions 
from Lemma \ref{explicit} yields a cubic polynomial in $a$; the order zero term clearly vanishes, while at order one, two and three it is straightforward to verify the identities
$$
(x_k+x_{k+1})\Big( x_0^{-1}x_{2k}-P^{(0)}\Big)
+x_0^{-1}x_1x_{2k}F_{-2k}^{(1)}+x_{2k}^{-1}x_0^2F_{2k+1}^{(1)}-x_{2k}F_{-2k+1}^{(1)}
=0,
$$ 
$$
F_{3k}^{(1)}+F_{3k+1}^{(1)}
+x_0^{-1}x_1x_{2k}F_{-2k}^{(2)}+F_{2k+1}^{(1)}F_{-2k}^{(1)}-x_{2k}F_{-2k+1}^{(2)}
-(x_k+x_{k+1})P^{(1)}
=0,
$$
$$ 
F_{3k+1}^{(2)}+(x_k+x_{k+1})\Big( x_0^{-1}F_{-2k}^{(2)}-P^{(2)}\Big)=0.
$$
This shows that $\varphi^*K=K$, and completes the proof of Theorem \ref{main}.
\end{prf}

\section{General solution and other corollaries}

\setcounter{equation}{0}

The result of Theorem \ref{main} has many implications for the recurrence (\ref{heiho}). 
Theorem \ref{mainhh} is just a special case: one finds the value $K=2k^2+8k+4$ by substituting 
$a=x_0=x_1=\ldots = x_{2k}=1$ into (\ref{kinv}). Also, note that we have not made use of the Laurent property in the proof, yet by the  formulae  in Lemma \ref{explicit} we see that $x_{-2k},x_{-2k+1},\ldots, x_{4k}$ provide $6k$ initial data for the linear recurrence ${\cal L}\, x_n =0$, and these are all Laurent polynomials in 
$x_0,x_1,\ldots,x_{2k}$ with coefficients in $\Z [a]$, as is $K$ given by  (\ref{kinv}), so we arrive at the following.
\begin{cor}\label{lp} The nonlinear recurrence (\ref{heiho}) has the Laurent property, i.e. 
$x_n \in \Z [x_0^{\pm 1}, x_1^{\pm 1}, \ldots, x_{2k}^{\pm 1},a]$ $\forall n\in\Z$. 
\end{cor} 

In addition to the homogeneous linear recurrences in Theorem \ref{lnrs}, various inhomogeneous linear recurrences now follow. 
\begin{cor}\label{inlins} 
The iterates of (\ref{heiho}) satisfy the linear recurrences 
\beq\label{in2k} 
x_{n+4k}-(K-1)\, x_{n+2k}+x_n = \nu_n, 
\eeq 
where $\nu_n$ is a $2k$-invariant, 
\beq \label{long} 
\sum_{j=0}^{2k-1} x_{n+4k+j}-(K-1)\, x_{n+2k+j}+x_{n+j}= K', 
\eeq 
where 
\beq\label{symfn} K' = \nu_0+\nu_1+\ldots + \nu_{2k-1} \eeq  
is a conserved quantity, and 
\beq\label{short} 
x_{n+2}+\eta_n \, x_{n+1}+\zeta_n \, x_n =\epsilon_n,
\eeq 
where $\epsilon_n, \zeta_n$ and $\eta_n$ are all $2k$-invariants. 
\end{cor} 
\begin{prfco} The operator $\cal L$ can be factorized as 
$$ 
{\cal L}=(\SH^{2k}-1)\Big(\SH^{4k}-(K-1)\SH^{2k}+1\Big) 
= (\SH-1)\left(\sum_{j=0}^{2k-1}\SH^j\right)\Big(\SH^{4k}-(K-1)\SH^{2k}+1\Big) ,
$$
which means that (\ref{gklin}) can be ``integrated'' in two different ways to yield (\ref{in2k}) and (\ref{long}). 
By shifting and summing $2k$ copies of (\ref{in2k}), the conserved quantity $K'$ is given as a symmetric function of the $2k$ independent shifts of $\nu_n$ according to (\ref{symfn}). 

The equation (\ref{gklin}) also implies that, for all $n$,  the vector $(1,-K,K,-1)^T$ belongs to the kernel of the matrix 
$$ 
\left(\begin{array}{cccc} 
1 & 1 & 1& 1 \\
 x_n & x_{n+2k} & x_{n+4k} & x_{n+6k} \\
x_{n+1} & x_{n+2k+1} & x_{n+4k+1}& x_{n+6k+1} \\
x_{n+2} & x_{n+2k+2} & x_{n+4k+2} & x_{n+6k+2}
\end{array}\right) 
$$ 
so by writing a vector in the kernel of its transpose as the row vector $(\epsilon_n,-\zeta_n,\eta_n,-1)$, the relation 
(\ref{short}) follows, and the same argument as in the proof of Theorem \ref{lnrs} shows that 
$\epsilon_n,\zeta_n,\eta_n$ are invariant under shifting $n\to n+2k$.
\end{prfco}

The fact that the iterates of (\ref{heiho}) satisfy a linear relation with constant coefficients means that they can be 
written explicitly in terms of the roots of the associated characteristic polynomial. Due to the particular form of (\ref{gklin}), 
the general solution can also be written  using either trigonometric functions or Chebyshev polynomials (with 
the latter form of the solution for $k=1$ being included in the results of \cite{honechapt}). 
In order to do this, we introduce quantities  $\theta$ and $t$  such that 
$$ 
t=\frac{K-1}{2}=\cos\Theta, \qquad \Theta = 2k\theta, 
$$ 
and recall that the Chebyshev polynomials of the first and second kinds are defined by 
$$
T_n(t)=\cos (n\Theta), \qquad U_{n-1}(t)=\frac{\sin(n\Theta)}{\sin\Theta}
$$ 
respectively, so that $T_0=U_0=1$, $T_1=T_{-1}=t$, $U_1=2t$, $U_{-1}=0$.
\begin{cor}
The general solution of (\ref{heiho}) can be written in the form 
$$ 
x_n =a_n +b_n \cos (n\theta +\phi_n), 
$$ 
where $a_n,b_n,\phi_n$ are all periodic in $n$ with periodic $2k$, or 
 as 
$$ 
x_n = q_j+r_j\,T_m(t)+s_j\,U_m(t), \qquad m =  \left \lfloor{\frac{n}{2k}}\right \rfloor , 
$$
where  $j=n\mod 2k$ and for $j=0,1,\ldots, 2k-1$ the coefficients are 
$$
\left(\begin{array}{c} q_j\\ r_j\\ s_j \end{array}\right) 
=\frac{1}{2t(1-t)} \, 
\left(\begin{array}{ccc} t &-2t^2 & t \\ -1 & 2t & 1-2t \\ 1-t & 0 & t-1 \end{array}\right) 
\left(\begin{array}{c} x_{2k+j}\\ x_j\\ x_{-2k+j} \end{array}\right). % , \quad j=0,1,\ldots, 2k-1. 
$$ 
\end{cor}

\section{Conclusions} 

\setcounter{equation}{0}

We have proved that all sequences generated by the nonlinear recurrence (\ref{heiho}) satisfy a linear relation, which 
was left as an open problem in \cite{heiho}. A key feature in the proof was to use the reversibility property 
in Lemma \ref{rev}. 
In fact, all the cluster maps obtained from period 1 quivers in \cite{fordymarsh09}, and many of the recurrences considered 
in \cite{alman}, are also reversible 
with the same sort of reversing symmetry, which means that the above approach can be applied 
in those cases too, and may prove useful for finding explicit formulae for conserved quantities (when they exist). 

One question that remains open is whether there is any natural 
Poisson structure which is preserved by the map (\ref{phi}), since Poisson structures (or presymplectic structures) arise 
naturally in the context of cluster algebras, but whether there is something similar for 
LP algebras in general  is an open question. In fact we expect that  there is a Poisson bracket (albeit a rather degenerate one, of rank two) for a combination of  two reasons: first, the map should have many conserved quantities in addition to $K$ and 
$K'$ given by (\ref{symfn}), since any symmetric function of the shifts of a $2k$-invariant is conserved; and second, 
$\varphi$ has the logarithmic volume form 
$$
\om = \frac{1}{x_0x_1\ldots x_{2k}}\, \rd x_0\wedge \rd x_1\wedge \ldots \wedge \rd x_{2k}
$$ 
which is anti-invariant, in the sense that $\varphi^*\om=-\om$; so 
if there are at least $2k-2$ independent conserved quantities, then the corresponding co-volume can be contracted with 
their differentials to construct a Poisson bracket, by one of the results in \cite{bhq}.

\noindent {\bf Acknowledgments.}  Some of these results first appeared in the PhD thesis \cite{ward}, which was supported by EPSRC studentship EP/P50421X/1. ANWH is supported by EPSRC fellowship EP/M004333/1.

\end{document}